\documentstyle[12pt,epsf,amssymb]{article}

\font\tenmsb=msbm10
\font\sevenmsb=msbm7
\font\fivemsb=msbm5
 
\newfam\msbfam 
\textfont\msbfam=\tenmsb
\scriptfont\msbfam=\sevenmsb
\scriptscriptfont\msbfam=\fivemsb
 
\def\Bbb#1{\fam\msbfam\relax#1}

\newtheorem{thm}{Theorem}[section]
\newtheorem{prop}[thm]{Proposition}
\newtheorem{cor}[thm]{Corollary}
\newtheorem{lem}[thm]{Lemma}
\newtheorem{conj}[thm]{Conjecture}
\newtheorem{exa}[thm]{Example}
\newtheorem{defn}[thm]{Definition}

\newtheorem{rem}[thm]{Remark}
\newtheorem{note}[thm]{Notation}
\newtheorem{alg}[thm]{Algorithm}

\newcommand{\ben}{\begin{enumerate}}
\newcommand{\een}{\end{enumerate}}
\newcommand{\ble}{\begin{lem}}
\newcommand{\ele}{\end{lem}}
\newcommand{\bth}{\begin{thm}}
\renewcommand{\eth}{\end{thm}}
\newcommand{\bpr}{\begin{prop}}
\newcommand{\epr}{\end{prop}}
\newcommand{\bco}{\begin{cor}}
\newcommand{\eco}{\end{cor}}
\newcommand{\bcon}{\begin{conj}}
\newcommand{\econ}{\end{conj}}
\newcommand{\bde}{\begin{defn}}
\newcommand{\ede}{\end{defn}}
\newcommand{\bex}{\begin{exa}}
\newcommand{\eex}{\end{exa}}
\newcommand{\brem}{\begin{rem}}
\newcommand{\erem}{\end{rem}}
\newcommand{\bnot}{\begin{note}}
\newcommand{\enot}{\end{note}}
\newcommand{\balg}{\begin{alg}}
\newcommand{\ealg}{\end{alg}}

\newcommand{\bib}{thebibliography}
\newcommand{\qed}{\square}

\newcommand{\C}{{\Bbb C}}

\newcommand{\PP}{{\Bbb P}}
\newcommand{\R}{{\Bbb R}}

\begin{document}

\title{A New Algorithm for Solving the Word Problem in Braid Groups} 

\author{D. Garber, S. Kaplan, M. Teicher$^1$ \\
\small Department of Mathematics and Computer Sciences\\
\small Bar-Ilan University\\
\small Ramat-Gan, Israel 52900 \{ garber,kaplansh,teicher \}@macs.biu.ac.il}

\stepcounter{footnote}\footnotetext{Partially supported by the Emmy Noether Research Institute for
Mathematics, Bar-Ilan University and the Minerva Foundation, Germany, and by the Excellency Center ``Group theoretic methods in
the study of algebraic varieties'' of the National Science Foundation of
Israel.}

\date{\today \\[1in]}
\maketitle

\begin{center}
proposed running head: Algorithm for Braid Word Problem

Author for proof: Shmuel Kaplan, 

Department of Mathematics and Computer Sciences

Bar Ilan University 

Ramat-Gan, Israel

Zip 52900
\end{center}

\begin{abstract}
One of the most interesting questions about a group is if its word problem can be solved and how. 
The word problem in the braid group is of particular interest to topologists, algebraists and geometers, 
and is the target of intensive current research.
We look at the braid group from a topological point of view (rather than a geometrical one). 
The braid group is defined by the action of
diffeomorphisms on the fundamental group of a punctured disk. We exploit the topological definition of 
the braid group in order to give a new approach for solving its word problem. 
Our algorithm is faster, in comparison with known algorithms, for short braid words with respect 
to the number of generators combining the braid, and it is almost independent of the number of strings in the braids. 
Moreover, the algorithm is based on a new computer presentation of the
elements of the fundamental group of a punctured disk. This presentation can be used also for other algorithms.
\end{abstract}

\begin{flushleft}
\small{Key Words: Fundamental group, Braid group, Word problem, Algorithm} 
\end{flushleft}

\tableofcontents

\section*{Introduction}
Let $D$ be a closed disk, and $K=\{k_1,...,k_n\}$ be $n$ points in $D$. Let $B$ be 
the group of all diffeomorphisms $\beta$ of $D$ such that $\beta(K)=K$, 
$\beta |_{\partial D}={\rm Id} |_{\partial D}$. The braid group is derived from $B$ by identifying two elements if their 
actions on $\pi _1(D \setminus K,u)$ are equal. To simplify the algorithm, we choose a
 geometric base of $\pi _1(D \setminus K,u)$, and
we look at the action of $\beta \in B$ on the elements of this geometrical base.

Thus, in order to determine if two words in the braid group are identical, we check 
whether their actions on the different elements of the chosen geometrical base are identical. 
Accordingly, to make this checking procedure efficient, we produced a new computerized presentation, and two new algorithms:

\ben
\item
A presentation of the geometrical base of $\pi _1(D \setminus K,u)$.
\item
An algorithm to compute the action.
\item
An algorithm for reducing the presentation into a unique form.
\een
The composition of these components holds the solution for the word problem in the braid group.

In section 1 we will give a short presentation of the fundamental group, algebraic 
and topological definitions of the braid group, and finally we will present the word problem 
in the braid group and some of the known solutions for it (Garside \cite{GAR}, Dehornoy \cite{Deh1}, 
Birman Ko and Lee \cite{NEW}). 
In section 2, we will present our algorithms. 
Section 3 will be dedicated to the proof of the correctness of the algorithms.
Section 4 will deal with some aspects of their complexity. 
Finally, in section 5 we give conclusions, future applications of the
new presentation, and further plans.

\section{Braid group and preliminaries}
In this section, we will recall some definitions that we will use in the sequel. 
Some of them will concern the fundamental group, others will
describe Artin's braid group. We will give two equivalent definitions of the braid group. 
The first definition is Artin's definition \cite{Artin}, and the second 
is based on the group of diffeomorphisms of the punctured disk. The latter
will give us the tools needed for solving the word problem, which will be
presented at the end of this section.

\subsection{The fundamental group}
$D$ is a closed oriented unit disk in $\R ^2$, $K=\{k_1,...,k_n\} \subset D$ 
is a finite set of points, and $u \in \partial D$. We look at the fundamental group of $D \setminus K$
denoted by $\pi _1(D \setminus K,u)$.
It is known that the fundamental group of a punctured disk with $n$ holes is a free group on $n$ generators.

Let $q$ be a simple path connecting $u$ with one of the $k_i$, say $k_{i_0}$, such that
$q$ does not meet any other point $k_j$ where $j \neq i_0$. To $q$ we will assign a
loop $l(q)$ (which is an element of $\pi _1(D \setminus K,u)$) as follows:

\bde{\underline{$l(q)$}}

Let $c$ be a simple loop equal to the (oriented) boundary of a small neighborhood $V$ of
$k_{i_0}$ chosen such that $q'=q \setminus (V \cap q)$ is a simple path.
Then \underline{$l(q)=q' \cup c \cup q'^{-1}$}. We will use the same notation for the element of $\pi
_1(D \setminus K,u)$ corresponding to $l(q)$.
\ede

\bde
Let  $(T_1,...,T_n)$ be an ordered set of simple paths in $D$ which connect the $k_i$'s with
$u$ such that:
\ben
\item
$T_i \cap k_j=\emptyset$ if $i \neq j$ for all $i,j=1,...,n$.
\item
$\displaystyle \bigcap_{i=1}^nT_i=\{u\}$.
\item
for a small circle $c(u)$ around $u$, each $u'_i=T_i \cap c(u)$ is a single point and the
order in $(u'_1,...,u'_n)$ is consistent with the positive orientation of $c(u)$.
\een

We say that two such sets $(T_1,...,T_n)$ and $(T'_1,...,T'_n)$ are \underline{equivalent} if 
$l(T_i)=l(T'_i)$ in $\pi _1(D \setminus K,u)$  for all $i=1,...,n$.

An equivalence class of such sets is called a \underline{bush} in $D \setminus K$.
\ede

\bde
A \underline{g-base} (geometrical base) of $\pi _1(D \setminus K,u)$ is an ordered 
free base of $\pi _1(D \setminus K,u)$
which has the form $(l(T_1),...,l(T_n))$, where $(T_1,...,T_n)$ is a bush in $D \setminus
K$.
\ede

For convenience, we choose $D$ to be the unit disk and the set $\{k_1,...,k_n\}$ on the $x$-axis 
ordered from left to right and $u=(0,-1)$ and hence $u \in \partial D$.

We would like to point out a particular g-base which will be used in the paper.
Choose $T_i$ to be the straight line connecting $u$ with $k_i$, then we call $(l(T_1),...,l(T_n))$ the 
\underline{standard g-base of $\pi _1(D \setminus K,u)$} and it is shown in the following figure: 

\begin{center}
\epsfysize=3cm
\epsfbox{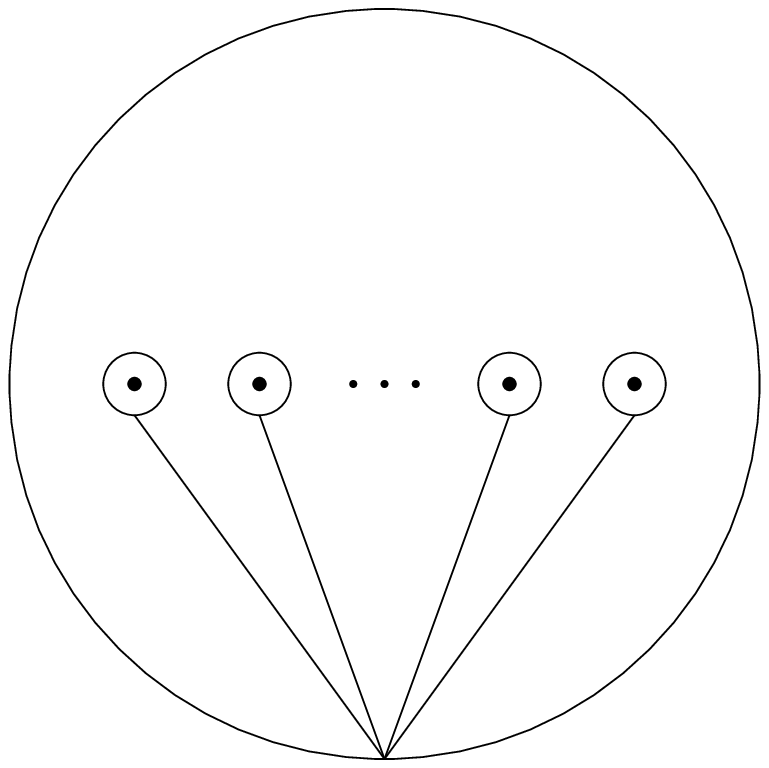}

\small{Figure $1$: The standard g-base}
\end{center}

\subsection{Artin's braid group}
In this subsection, we will give two equivalent definitions for the braid group. The first
is algebraic and the second is topological, which will be used to present our algorithms in this paper.

\subsubsection{The algebraic definition for the braid group}

Here we will lay out Artin's definition \cite{Artin} as used in most cases.

\bde
\underline{Artin's braid group $B_n$} is the group generated by $\{\sigma _1,...,\sigma _{n-1}\}$
submitted to the relations 
\ben
\item $\sigma _i\sigma _j=\sigma _j\sigma _i$ where $|i-j| \geq 2$ 
\item $\sigma _i\sigma _{i+1}\sigma _i=\sigma _{i+1}\sigma _i\sigma _{i+1}$  for all $i=1,...,n-2$
\een
\ede

One can look at this as a geometrical definition, since it can be interpreted to the set of 
ties of $n$ strings going from top to bottom. This is done by assigning a positive 
switch between any adjacent pair of strings to one of the generators. This means that $\sigma _i$ corresponds 
to the geometrical element described in the following figure:

\begin{center}
\epsfxsize=2cm
\epsfbox{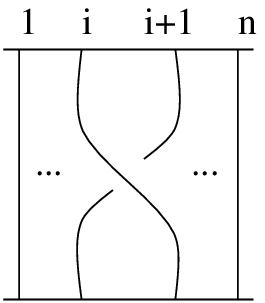}

\small{Figure$2$: The geometrical braid associated with $\sigma _i$}
\end{center}

The operation for this group can be described as the concatenation of two geometrical 
sets of strings resulting in what is called a \underline{braid}.

\bex
The geometrical braid that corresponds to $\sigma _1 \sigma _2^{-1} \sigma _1 \sigma _3$ 
is presented in the following figure:

\begin{center}
\epsfxsize=2cm
\epsfbox{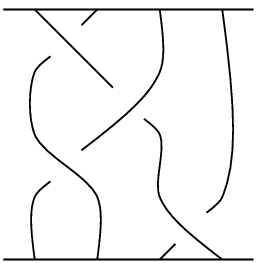}

\small{Figure $3$: The geometrical braid $\sigma _1 \sigma _2^{-1} \sigma _1 \sigma _3$}
\end{center}
\eex

\subsubsection{The topological definition for the braid group}

Let $D,K,u$ be as above.

\bde
Let $B$ be the group of all diffeomorphisms $\beta$ of $D$ such that $\beta(K)=K$, $\beta
|_{\partial D}={\rm Id}|_{\partial D}$. For $\beta _1,\beta _2 \in B$ we say that $\beta _1$ is
\underline{equivalent} to $\beta _2$ if $\beta _1$ and $\beta _2$ induce the same automorphism of $\pi
_1 (D \setminus K,u)$. The quotient of $B$ by this equivalence relation is called 
\underline{the braid group $B_n[D,K]$} ($n=\#K$). The elements of $B_n[D,K]$ are called \underline{braids}.
\ede

\brem
For the canonical homomorphism $\psi:B \to Aut(\pi _1(D \setminus K,u))$,
we actually have $B_n[D,K] \cong Im(\psi)$. 
\erem

We recall two facts from \cite{BGTI}[section III].

\ben
\item
If $K' \subset D'$, where $D$ is another disk, and $\#K=\#K'$ then $B_n[D,K] \cong B_n[D',K']$.
\item
Any braid $\beta \in B_n[D,K]$ transforms a g-base to a g-base. Moreover, 
for every two g-bases, there exists a unique braid which transforms one g-base to another.
\een

We distinguish some elements in $B_n[D,K]$ called \underline{half-twists}.

Let $a,b \in K$ be two points. We denote $K_{a,b}=K \setminus
\{a,b\}$. Let $\sigma$ be a simple path in $D \setminus (\partial D \cup K_{a,b})$ connecting
$a$ with $b$. Choose a small regular neighborhood $U$ of $\sigma$ and an orientation
preserving diffeomorphism $f:\R ^2 \to \C$ such that $f(\sigma )=[-1,1]$, $f(U)=\{z \in \C \ | \ |z| <2\}$. 

Let $\alpha (x)$, $0 \leq x$ be a real smooth monotone function such that:

$$\alpha (x)=\left \{ \matrix{1 & 0 \leq x \leq \frac{3}{2} \cr
                             0 & 2 \leq x} \right.$$

Define a diffeomorphism $h:\C \to \C$ as follows: for $z=re^{i\varphi} \in \C$ let
$h(z)=re^{i(\varphi +\alpha (r)\pi )}$

For the set $\{z \in \C \ | \ 2 \leq |z|\}$,  $h(z)={\rm Id}$,
and for the set $\{z \in \C \ | \ |z|\leq \frac{3}{2}\}$, $h(z)$ a rotation by 
$180 ^{\circ}$ in the positive direction.

The diffeomorphism $h$ defined above induces an automorphism on $\pi _1(D \setminus K,u)$,
that switches the position of two generators of $\pi _1(D \setminus K,u)$, as
can be seen in the figure: 

\begin{center}
\epsfysize=3cm
\epsfbox{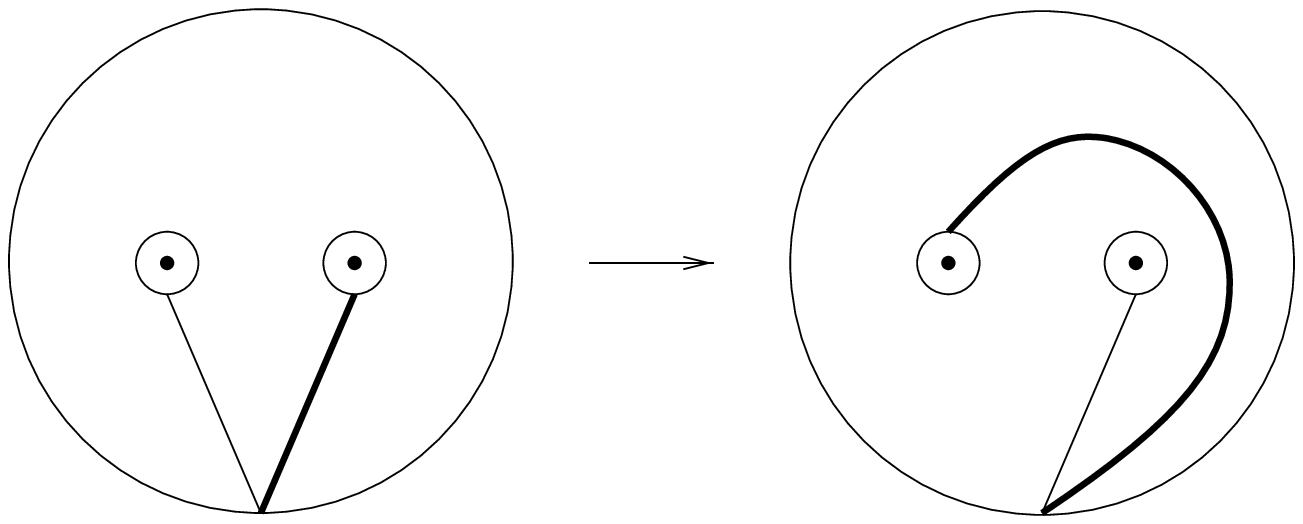}

\small{Figure $4$: The action of the diffeomorphism $h$}
\end{center}

Considering $(f \circ h \circ f^{-1})|_D$ (we will compose from left to
right), we get a diffeomorphism of $D$ which switches $a$ and $b$ and is the identity on $D
\setminus U$. Thus it defines an element of $B_n[D,K]$.

\bde
Let $H(\sigma)$ be the braid defined by $(f \circ h \circ f^{-1})|_D$. We call
$H(\sigma )$ the \underline{positive half-twist defined by $\sigma$}.
\ede

The half-twists generate $B_n$. In fact, one can choose $n-1$ half-twists that generates $B_n$ (see below):

\bde
Let $K=\{k_1,...,k_n\}$, and $\sigma _1,...,\sigma _{n-1}$ be a system of simple paths in
$D \setminus \partial D$ such that each $\sigma _i$ connects $k_i$ with $k_{i+1}$ and 

for all $i,j \in \{1,...,n-1\}$, $i<j$, $\left \{ 
\matrix{\sigma _i \cap \sigma _j=\emptyset & 2 \leq |i-j| \cr
\sigma _i \cap \sigma _{i+1}=\{k_{i+1}\} & i=1,...,n-2} \right.$. 

Let $H_i=H(\sigma _i)$. The ordered system of
(positive) half twists $(H_1,...,H_{n-1})$ are called a \underline{frame of $B_n[D,K]$}.
\ede

\bth
If $(H_1,...,H_{n-1})$ is a frame of $B_n[D,K]$, then
$B_n[D,K]$ is generated by $\{H_i\}_{i=1}^{n-1}$. Moreover, 
if $(H_1,...,H_{n-1})$ is a frame of $B_n[D,K]$, then the set
$\{H_i\}_{i=1}^{n-1}$ with the two relations $H_iH_j=H_jH_i$ if $2 \leq |i-j|$ and
$H_iH_{i+1}H_i=H_{i+1}H_iH_{i+1}$ for any $i=1,...,n-2$ are sufficient enough to present 
$B_n[D,K]$ and therefore this definition and Artin's definition for the braid group are equivalent.
\eth

\noindent
{\it Proof:} See \cite{BGTI}.

As the \underline{standard frame} we will use a frame which its paths are the straight segments connecting the 
point $k_i$ to $k_{i+1}$ $i=1,...,n-1$.

\subsubsection{The word problem}

First we define what is called a braid word.

\bde
Let $b \in B_n$ be a braid. Then it is clear that $b=\sigma _{i_1}^{e_1} \cdot ...\cdot \sigma
_{i_l}^{e_l}$ for some sequence of generators, where $i_1,...,i_l \in \{1,...,n-1\}$ and 
$e_1,...,e_l \in \{1,-1\}$. We will call such a presentation of $b$ a \underline{braid
word}, and $\sigma _{i_k}^{e_k}$ will be called the \underline{$k^{th}$ letter of the word $b$}. 
$l$ is the \underline{length} of the braid word.
\ede

We will distinguish between two relations on the braid words.

\bde
Let $w_1$ and $w_2$ be two braid words. We will say that $w_1=w_2$ if they represent the
same element of the braid group.
\ede

\bde
Let $w_1$ and $w_2$ be two braid words. We will say that $w_1 \equiv w_2$ if $w_1$ and
$w_2$ are identical letter by letter.\
\ede

Now, we can introduce the word problem: Given two braid words $w_1$ and $w_2$, 
decide whether $w_1=w_2$ or not.

\subsection{Two known algorithms for the word problem}
There are several known algorithms for solving the word problem for the braid group. 
In this section, we will summarize 
some of them. The complexity of different algorithms varies, 
but to our knowledge, the best known solution is of complexity of $O(l^2)$, 
where $l$ is the length of the longer braid word.

\subsubsection{Garside's solution}
Garside \cite{GAR} gave a solution for the braid word problem in 1969. 
His solution is based on the definition of positive words, which contain only
generators with positive power. Then, he stated that the fundamental word of the braid group
$\Delta _n$ has a property that enables to replace all the generators with a negative power. 
This can be done simply by noticing the fact that for any $i$, there exists a
positive braid word $w_i$ for which $\sigma _i ^{-1}=\Delta _n^{-1}w_i$. 

Another property of the $\Delta _n$ is that for any $i$ we have
that $\sigma _i\Delta _n=\Delta _n\sigma _{n-i}$. This gives a method for writing a given
braid word $w$ in such a way that $w=w_1w_2$ where $w_1=\Delta _n^r$, $r \leq 0$ which is a negative 
braid word and $w_2$ is a positive braid word.

Now, one can write $w_2=\Delta _n^qw_3$, where $q$ is maximal. By doing this,
he can increase $r$ resulting in the minimal way of writing $w=\Delta
_n^{r-q}w_3$. By organizing $w_3$ in a lexicographic order, we obtain what is called \underline{Garside's
normal form} of the braid word $w$. 

Garside proved, that two braid words $w$ and $w'$ are equal if and only if their normal forms are the
same.

There are some implementations for solving the braid word using this solution, and variations of it 
as can be found, for example, in \cite{POSBR}, \cite{BAND}, \cite{NEW} and \cite{EFF}. 
For achieving the best complexity by this method, one has to expand the size of the set of 
generators of the braid group, resulting in the complexity of $O(l^2)$ where $l$ is the length of the
longer of the two braid words.

\subsubsection{Dehornoy's solution}
Dehornoy (\cite{Deh1}, \cite{Deh2}) used 
a different approach for solving the problem. His approach is based on a 
definition of a $\sigma $-reduced braid word, which is a braid word that for any integer $i$, any
occurrence of the letter $\sigma _i$ is separated from any occurrence of the letter $\sigma
_i^{-1}$ by at least one occurrence of a letter $\sigma _j^{\pm1}$ with $j<i$.

Dehornoy presented an algorithm for transforming any braid word to its reduced form. He
proved that the reduced form of a braid word $w$ is ${\rm Id}$ (i.e. the null braid word) 
if and only if $w$ is the identity word.
This gives a simple way of checking whether two braid words $w$ and $w'$ are equal, simply 
by writing $w''=w(w')^{-1}$ and reducing $w''$. If the reduced form of $w''$ is ${\rm Id}$, it means that $w=w'$. 

The reduction process is actually a type of an unknotting process that unties the twisted strings in 
a braid, by adding proper sequences and transforming locally twisted strings into an untwisted
state as shown in the following figure:

\begin{center}
\epsfysize=6cm
\epsfbox{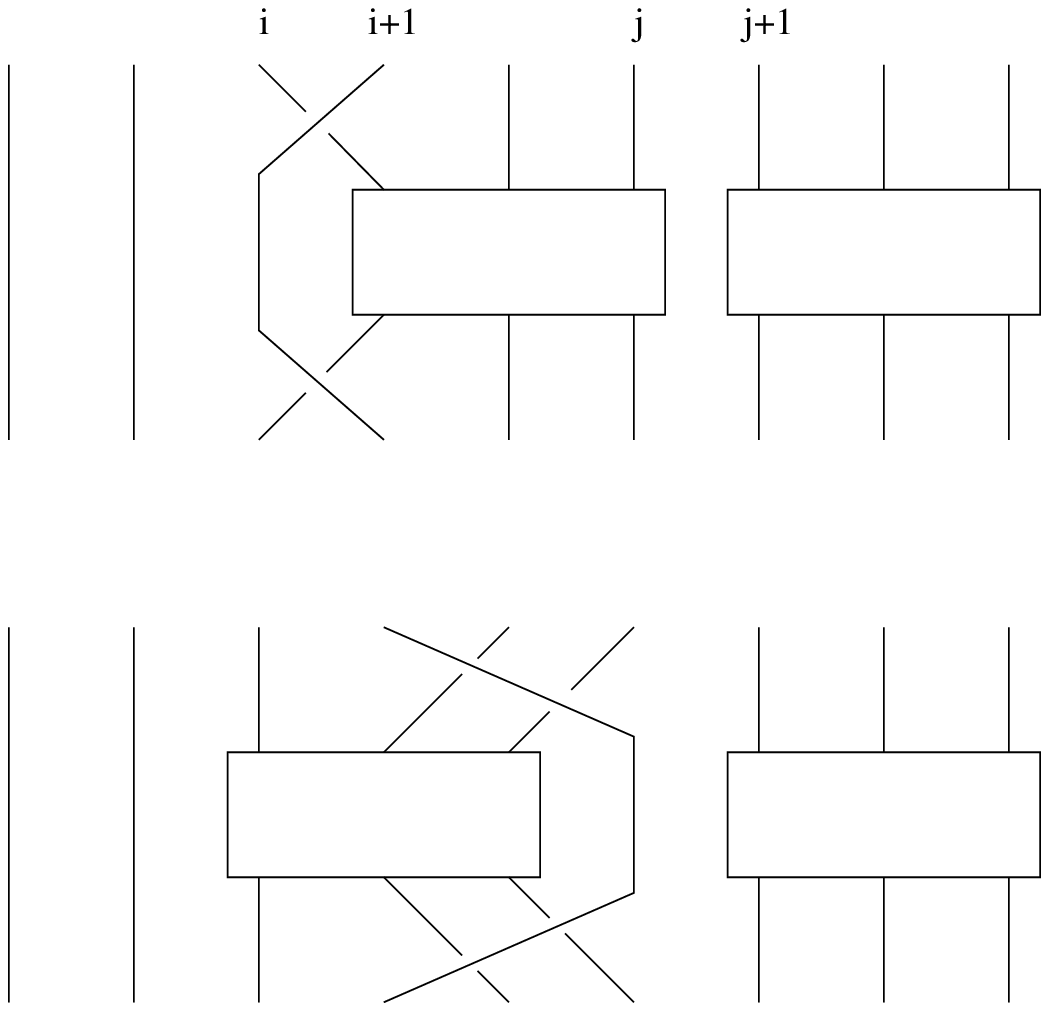}

\small{Figure $5$: Unknotting process in Dehornoy's algorithm}
\end{center}

Dehornoy conjectured that the complexity of his algorithm is bounded by $O(l^2)$ where $l$ is the length of
the longer braid word.

In the next section we will present our algorithm, which is based on a completely different
approach. 

\section{The presentation of the new algorithm for solving the word problem}
The algorithm that we are going to present in order to solve the word problem in the braid group 
is based on the interplay between its two definitions. 
We will fix the standard frame and the standard g-base that will be used as a starting position.  
We associate the generator $\sigma _i$ to the half-twist $H_i$ in the standard frame for every $i=1,...,n-1$.
By using our two algorithms and encoding the g-bases in a unique way, and by using an algorithmic way to explore the 
changes that happen to the standard g-base while the braid word acts on it, 
we produce a practical algorithm for the word problem. Mathematically, we compare two braid words by  
taking  one braid word and compute the result of its action on the standard g-base of the fundamental 
group. Then, we take the other braid word and compute the same result. 
The two braid words are equal if and only if the two resulted g-bases are identical.

\subsection{The computerized implementation of the g-base}
In this subsection, we will describe the way we encode the g-base. It involves some
conventions.

Recall that $D$ is the closed unit disk, the point $u$ is the point $(0,-1)$ and the points
in $K$ are on the $x$-axis.

In order to encode the path in $D$, which is an element of the g-base, we will
distinguish some positions in $D$.

\bnot
We will denote by $(i,1)$ a point {\bf close} to $k_i$ but above it, $(i,-1)$ a point {\bf close} to 
$k_i$ but below it, and $(i,0)$ the point $k_i$ itself.

We will also denote the point $u$ by $(-1,0)$ (which is not its position in $D$, rather 
only a notation).
\enot

To represent a path in $D$, we will use a linked list which its links are based on the notations above, which
represents the position of the path in relation to the points $u$ and $k_i$, $i=1,...,n$.

Each link of the list holds the two numbers as described above. We will call them
(point,position).

\bex
The list $(1,0) \to (2,1) \to (3,1) \to (4,-1) \to (5,0)$ represents the following path:
\eex

\begin{center}
\epsfysize=1cm
\epsfbox{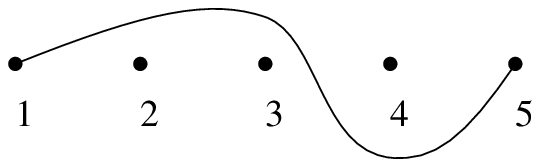}

\small{Figure $6$}
\end{center}

As a rule, we will never connect the point $u$ to any point $(i,-1)$. This will
be done in order to obtain a unique way of representation, and to make the automatic 
computation of the twists easier.

We will be able to tell whether a path $(-1,0) \to (i,1)$ is passing to the left or to the 
right of the point $i$ simply by checking its continuation. 
If the path is turning to the left $(-1,0) \to (i,1) \to (i-1,e)$, then it is passing to the right
of the point $i$, and if the path is turning to the right $(-1,0) \to (i,1) \to (i+1,e)$, then it
is passing to the left of the point $i$ (where $e \in \{-1,1,0\}$).

\bex
The list $(-1,0) \to (3,1) \to (2,0)$ represents the following path:

\begin{center}
\epsfysize=3cm
\epsfbox{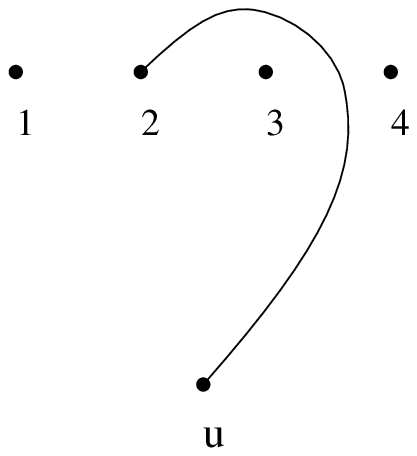}

\small{Figure $7$}
\end{center}

The list $(-1,0) \to (2,1) \to (3,1) \to (3,-1) \to (2,0)$ represents the following path:

\begin{center}
\epsfysize=3cm
\epsfbox{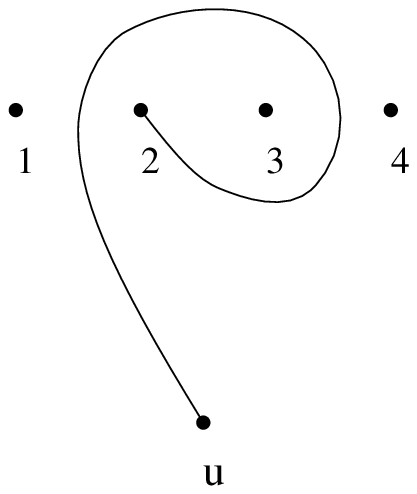}

\small{Figure $8$}
\end{center}
\eex

In order to unify our treatment of all the paths of the g-base, we will concatenate all 
of them into one list, which means that after we arrive at the end of one path (i.e. a
link $(i,0)$), the following link will be $(-1,0)$ marking the beginning of the next path.
For convenience, and not for mathematical reasons, we add the link $(-1,0)$ at the end of the
list.

\bex
The list $(-1,0) \to (1,1) \to (2,0) \to (-1,0) \to (1,0) \to (-1,0) \to (4,0) \to (-1,0) 
\to (4,1) \to (3,0) \to (-1,0)$ represents the g-base in the following figure 
(the small circles around the points are omitted):

\begin{center}
\epsfysize=2.5cm
\epsfbox{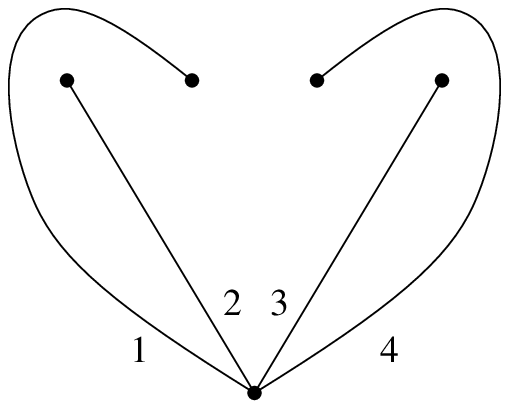}

\small{Figure $9$: The g-base represented by the list}
\end{center}
\eex

\subsection{The algorithm}
Now, we are ready to present the algorithm:

\balg
\parindent0pt
\parskip0pt

ProcessWord($w$)

{\bf input}: $w$ - a braid word.

{\bf output}: a list which represents the g-base resulted after the action of the word's
letters on the standard g-base.

{\bf ProcessWord($w$)}

$g \leftarrow$ list that represents the standard g-base.

For every letter $\sigma _i$ in $w$ do
\ben
\item
act on $g$ using $\sigma _i$ by applying $(Positive/Negative)HalfTwist(\sigma _i,g)$
function.
\item
reduce $g$ to its unique form using $Reduce(g)$ function.
\een

return $g$.
\ealg
\parindent10pt
\parskip10pt

Now, we will present the $PositiveHalfTwist(\sigma _i,g)$ function.

\balg
\parindent0pt
\parskip0pt

PositiveHalfTwist($\sigma _i,g$)

{\bf input}: 

$\sigma _i$ - the generator of the braid group acting on the g-base.

$g$ - the list representing the g-base.

{\bf output}: a list representing the g-base after the action of $\sigma _i$ on $g$

{\bf PositiveHalfTwist($\sigma _i,g$)}

for each sequence of links in $g$ of the type $(i,e)$ or $(i+1,e)$ $(e \in \{-1,1,0\})$ do

$\qquad$ $BeforeSection \leftarrow$ The link just before the first link in the sequence

$\qquad$ $AfterSection \leftarrow$ The link after the last link in the sequence

$\qquad$ $FirstLink \leftarrow$ The first link in the sequence

$\qquad$ $SecondLink \leftarrow$ the second link in the sequence

$\qquad$ if $BeforeSection=(-1,0)$ then

$\qquad$ $\qquad$ act upon one of the following cases:

$\qquad$ $\qquad$ if $FirstLink.Point=i$ and $SecondLink.Point=0$ then

$\qquad$ $\qquad$ $\qquad$ add the link $(i-1,-1)$ after $BeforeSection$

$\qquad$ $\qquad$ $\qquad$ $BeforeSection \leftarrow$ the new link

$\qquad$ $\qquad$ if $FirstLink.Point=i+1$ and $FirstLink.Position=0$ then

$\qquad$ $\qquad$ $\qquad$ add the link $(i+2,-1)$ after $BeforeSection$

$\qquad$ $\qquad$ $\qquad$ $BeforeSection \leftarrow$ the new link

$\qquad$ $\qquad$ if $FirstLink.Point=i$ and $SecondLink.Point=i+1$ then

$\qquad$ $\qquad$ $\qquad$ add the link $(i-1,-1)$ after $BeforeSection$

$\qquad$ $\qquad$ $\qquad$ $BeforeSection \leftarrow$ the new link

$\qquad$ $\qquad$ if $FirstLink.Point=i$ and $SecondLink.Point=i-1$ then

$\qquad$ $\qquad$ $\qquad$ add the links $(i-1,-1) \to (i,-1)$ after $BeforeSection$

$\qquad$ $\qquad$ $\qquad$ $BeforeSection \leftarrow$ the first new link

$\qquad$ $\qquad$ if $FirstLink.Point=i+1$ and $SecondLink.Point=i+2$ then

$\qquad$ $\qquad$ $\qquad$ add the links $(i+2,-1) \to (i+1,-1)$ after $BeforeSection$

$\qquad$ $\qquad$ $\qquad$ $BeforeSection \leftarrow$ the first new link

$\qquad$ $\qquad$ if $FirstLink.Point=i+1$ and $SecondLink.Point=i$ then

$\qquad$ $\qquad$ $\qquad$ add the link $(i+2,-1)$ after $BeforeSection$

$\qquad$ $\qquad$ $\qquad$ $BeforeSection \leftarrow$ the new link

$\qquad$ for any link $L$ between $BeforeLink$ and $AfterLink$ do

$\qquad$ $\qquad$ $L.Position \leftarrow -L.Position$

$\qquad$ $\qquad$ $L.Point \leftarrow 2i+1-L.Point$

$\qquad$ if $BeforeSection.Point=i-1$ then

$\qquad$ $\qquad$ add the links $(i,-1) \to (i+1,-1)$ after $BeforeSection$

$\qquad$ else

$\qquad$ $\qquad$ add the links $(i+1,1) \to (i,1)$ after $BeforeSection$

$\qquad$ if $AfterSection.Point=i-1$ then

$\qquad$ $\qquad$ add the links $(i+1,-1) \to (i,-1)$ after $BeforeSection$

$\qquad$ else

$\qquad$ $\qquad$ add the links $(i,1) \to (i+1,1)$ after $BeforeSection$
\ealg
\parindent10pt
\parskip10pt

In order to obtain the $NegativeHalfTwist(\sigma _i,g)$ function, one has 
to use the $PositiveHalfTwist(\sigma _i,g)$ function while replacing the 
last two 'if statements' with the following:

\parindent0pt
\parskip0pt
{\it
$\qquad$ if $BeforeSection.Point=i-1$ then

$\qquad$ $\qquad$ add the links $(i,1) \to (i+1,1)$ after $BeforeSection$

$\qquad$ else

$\qquad$ $\qquad$ add the links $(i+1,-1) \to (i,-1)$ after $BeforeSection$

$\qquad$ if $AfterSection.Point=i-1$ then

$\qquad$ $\qquad$ add the links $(i+1,1) \to (i,1)$ after $BeforeSection$

$\qquad$ else

$\qquad$ $\qquad$ add the links $(i,-1) \to (i+1,-1)$ after $BeforeSection$
}
\parindent10pt
\parskip10pt

Now, we will present the algorithm for the function $Reduce(g)$. This function 
reduces the list that represents the g-base to a unique form without 
changing its homotopy type. This is done by applying several reduction rules 
that are induced from homotopic equivalences. The full proof of the validity 
of the rules will be given in the next section.

\balg
\parindent0pt
\parskip0pt

Reduce($g$)

{\bf input}: $g$ - a list representing a g-base.

{\bf output}: a list which represents a g-base homotopic to $g$. Its representation is unique.

{\bf Reduce($g$)}

for each link $L$ in the list do

$\qquad$ $FirstLink \leftarrow$ the link right after $L$

$\qquad$ $SecondLink \leftarrow$ the link right after $FirstLink$

$\qquad$ if $FirstLink=SecondLink$ then

$\qquad$ $\qquad$ delete $FirstLink$ and $SecondLink$ from the list

$\qquad$ if $FirstLink=(i,1)$ or $(i,-1)$ and $SecondLink=(i,0)$ then

$\qquad$ $\qquad$ delete FirstLink from the list

$\qquad$ if $FirstLink=(i,0)$ then

$\qquad$ $\qquad$ delete all links between $FirstLink$ and the first appearance 
of $(-1,0)$

$\qquad$ if $FirstLink=(-1,0)$ then

$\qquad$ $\qquad$ delete all links of the type $(i,-1)$ after it

$\qquad$ $L \leftarrow$ the next or previous link as necessary
\ealg
\parindent10pt
\parskip10pt

\section{Verification of the new algorithm (correctness)}
In this section, we will lay out the proof for the correctness of the two algorithms.

\subsection{Correctness of the $(Positive/Negative)HalfTwist(\sigma _i,g)$ algorithm}
We will begin our proof of the correctness of the algorithm by proving that 
the algorithm works on parts of the paths that are not directly connected to $u$ (i.e.
$(-1,0)$ is not BeforeSection).

\bpr
Let $\sigma _i$ be the generator acting on the g-base. Then any part of the
path which does not contain the points $i$ or $i+1$ is not affected by the twist.
\epr

\noindent
{\it Proof:} Since the action of the twist is defined locally, any part of the
path out of the twisted region (that contains only the points $i$ and $i+1$) is not
affected.
\hfill $\qed$

We need to check the behavior of the path locally in the twisted region. By 
\underline{local behavior} we mean the behavior of the links of the type $(i,e)$ or $(i+1,e)$, 
where $\sigma _i$ is the generator of the specified letter in the braid word, and 
$e\in\{-1,1,0\}$

\bpr
Let $\sigma _i$ be the generator acting on the g-base. The local behavior of the path is given by the 
following changes:
\ben
\item
The link's position changes to $-$position
\item
The link's point changes from $i$ to $i+1$ and vice versa.
\een
\epr 

\noindent
{\it Proof:} From its definition, the actual local action of the braid is a rotation of $180^\circ$. 
Therefore, a part of the path of the g-base's element, which was beneath a point 
before the rotation, will now be above a point, and the part of the path that 
was above a point before the rotation will now be beneath a point. Hence, if the position 
was equal to $-1$ before the twist, it will be equal to $1$ after the twist, and 
vice versa.

Moreover, if the point in the path was equal to $i$, then the point in the path will 
be $i+1$ after the twist, and if the point in the path was equal to $i+1$, then the point 
in the path will be $i$ after the twist.
\hfill $\qed$  

After we have rotated the path locally, we will have to connect it to the global path. 
This should be done by adding proper prefix and postfix sequences before and after the 
part that has been twisted.

\bpr \label{prefix}
Let $\sigma _i$ be the positive half-twist acting on the g-base. Then, 
the prefix sequence we have to add is as follows:
\ben
\item
$(i,-1) \to (i+1,-1)$ if the local section of the path is connected to a point to the 
left of the point $i$.
\item
$(i+1,1) \to (i,1)$ if the local section of the path is connected to a point to the right 
of the point $i+1$.
\een
\epr

\noindent
{\it Proof:} If the point just before the local section of the path is to the left of the 
twist, then the connecting path should be beneath the twisted region. On the contrary, if the 
point just before the local section is to the right of the twist, then the connecting path 
should be above the twisted region. So, all we need to add is the two links above 
the twisted region or beneath it as necessary, as shown in the following figure:

\begin{center}
\epsfxsize=12cm
\epsfbox{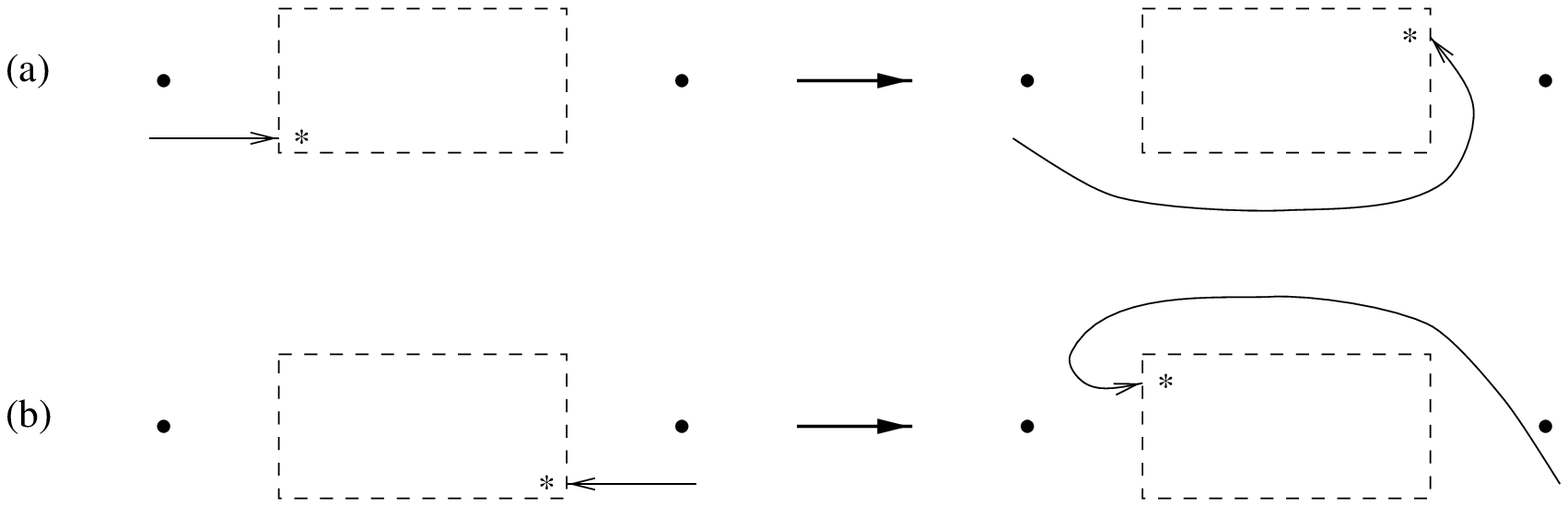}

\small{Figure $10$: Prefix added after the twist}
\end{center}

\hfill $\qed$

\bpr
Let $\sigma _i$ be the positive half-twist acting on the g-base. Then, the 
postfix sequence we have to add is as follows:
\ben
\item
$(i+1,-1) \to (i,-1)$ if the local section of the path is connected to a point left to 
the point $i$.
\item
$(i,1) \to (i+1,1)$ if the local section of the path is connected to a point to the right 
of the point $i+1$.
\een
\epr

\noindent
{\it Proof:} The proof is similar to the proof of proposition \ref{prefix}, see the following figure:

\begin{center}
\epsfxsize=12cm
\epsfbox{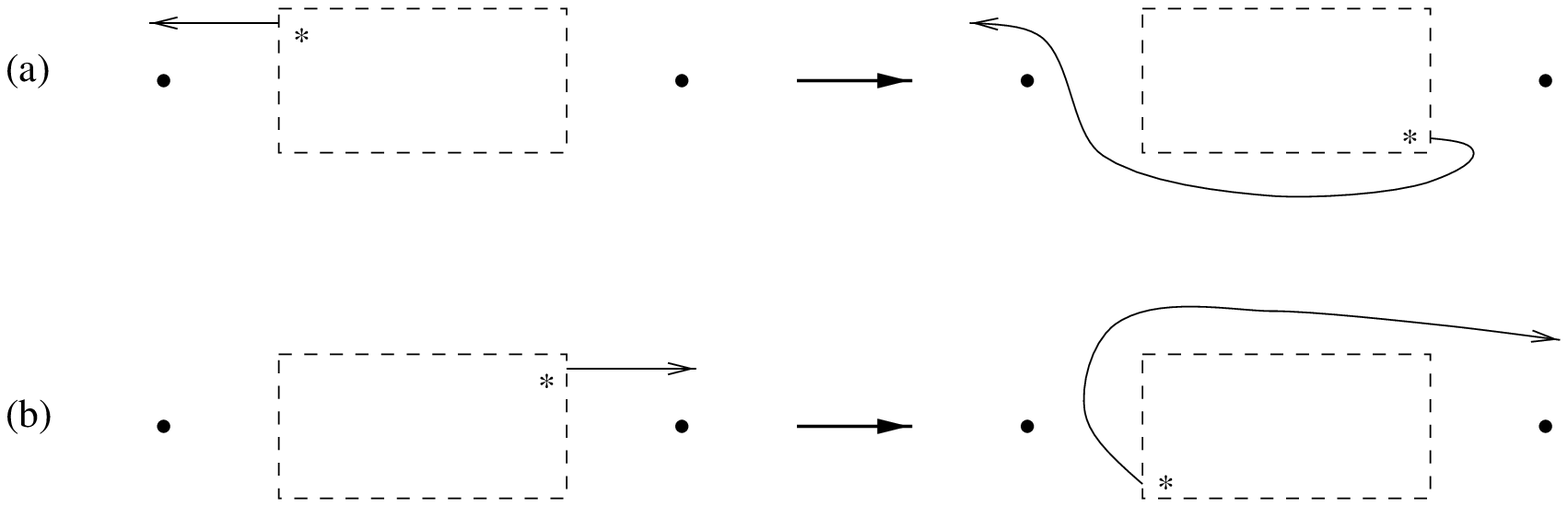}

\small{Figure $11$: Postfix added after the twist}
\end{center}

\hfill $\qed$

The local action of the braid generator $\sigma _i^{-1}$ is computable in the same way as the 
action of the generator $\sigma _i$. The prefix and the postfix sequences that we have to add 
are not the same sequences due to the direction of the twist.
Therefore, we have the following proposition:

\bpr
Let $\sigma _i$ be the negative half-twist acting on the g-base. Then, the 
prefix sequence we have to add is as follows:
\ben
\item
$(i,1) \to (i+1,1)$ if the local section of the path is connected to a point to the left of 
the point $i$.
\item
$(i+1,-1) \to (i,-1)$ if the local section of the path is connected to a point to the right 
of the point $i+1$.
\een
The postfix sequence we have to add is as follows:
\ben
\item
$(i+1,1) \to (i,1)$ if the local section of the path is connected to a point left of the point $i$.
\item
$(i,-1) \to (i+1,-1)$ if the local section of the path is connected to a point right of the point $i+1$.
\een
\epr

Now, we will consider the case where the link $(-1,0)$ is followed immediately by 
the local section of the path. In this case, we alter the path homotopically 
so that the preceding link to the local section of the path will not be $(-1,0)$. 
By doing this, we will reduce the problem to the one already proved by the above propositions, 
hence, we will be able to use the same algorithmic methods in these cases.

We have $6$ possible different cases:
\ben
\item
If we have the sequence $(-1,0) \to (i,0)$, then we add a link just below the point to 
the left of the local section which is $(i-1,-1)$. As a result, this point is the 
one preceding the local section (see figure (a)).
\item 
If we have the sequence $(-1,0) \to (i+1,0)$, then we add a link just below the point to the
right of the local section which is $(i+2,-1)$. As a result, this point is the one
preceding the local section (see figure (b)).
\item
If we have the sequence $(-1,0) \to (i,1) \to (i+1,e)$ $(e \in \{-1,1,0\})$, 
then we add a link just below the point to the left of the local section which is 
$(i-1,-1)$. As a result, this point is the one preceding the local section 
(see figure (c)).
\item
If we have the sequence $(-1,0) \to (i,1) \to (i-1,e)$ $(e \in \{-1,1,0\})$, 
then we add two links. The first is just below the point to the left of the local section 
which is $(i-1,-1)$ and therefore will be the preceding of the local section sequence, and
the second will be just below the point $i$ which is $(i,-1)$
(see figure (d)).
\item
If we have the sequence $(-1,0) \to (i+1,1) \to (i+2,e)$ $(e \in \{-1,1,0\})$, 
then we add two links. The first is just below the point to the right of the local section 
which is $(i+2,-1)$ and therefore will be the preceding of the local section sequence, and
the second will be just below the point $i+1$ which is $(i+1,-1)$
(see figure (e)).
\item
If we have the sequence $(-1,0) \to (i+1,1) \to (i,e)$ $(e \in \{-1,1,0\})$, 
then we add a link just below the point to the right of the local section which is 
$(i+2,-1)$. As a result, this point is the one preceding the local section 
(see figure (f)).
\een  

\begin{center}
\epsfysize=8cm
\epsfbox{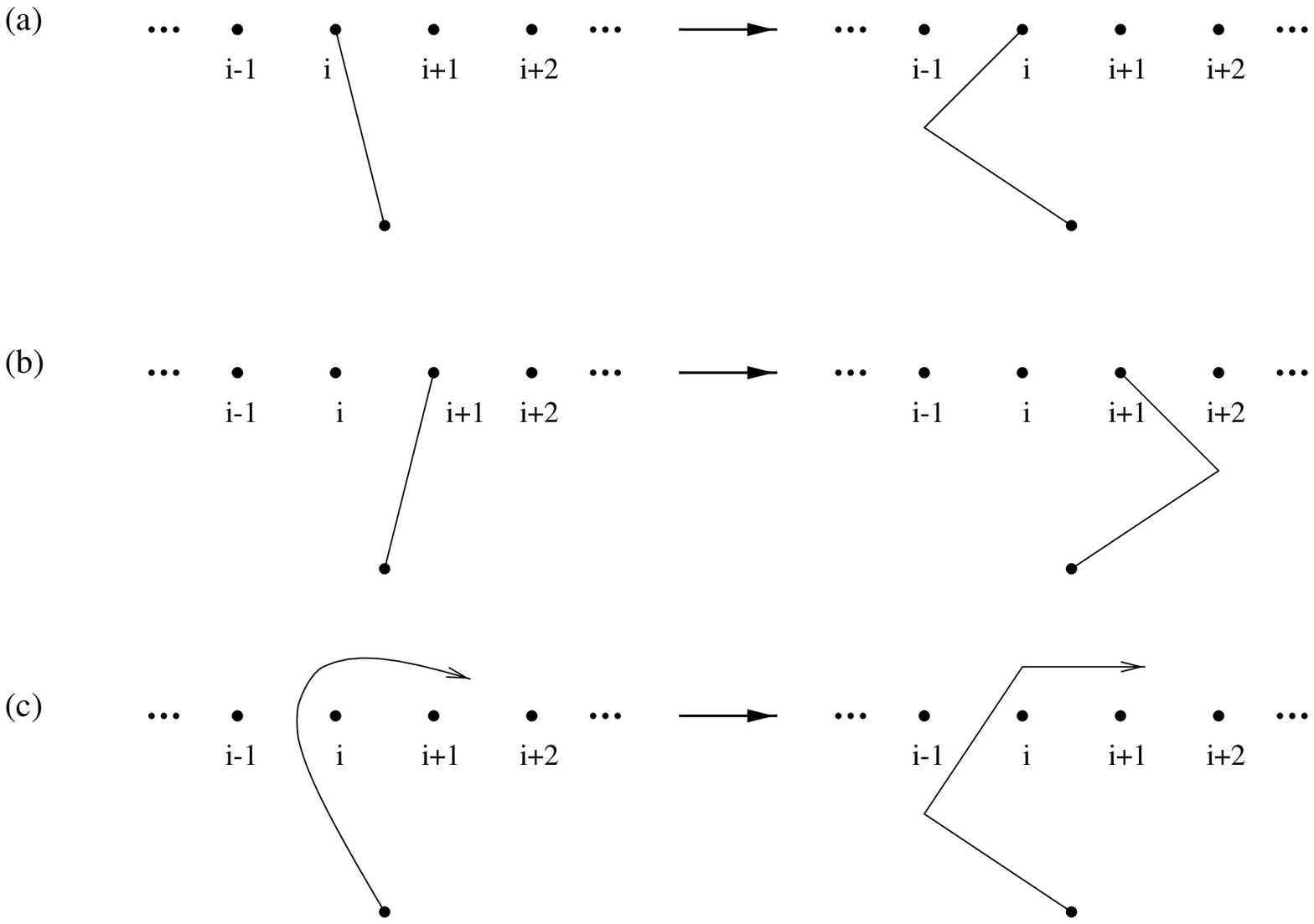}
\end{center}

\begin{center}
\epsfysize=8cm
\epsfbox{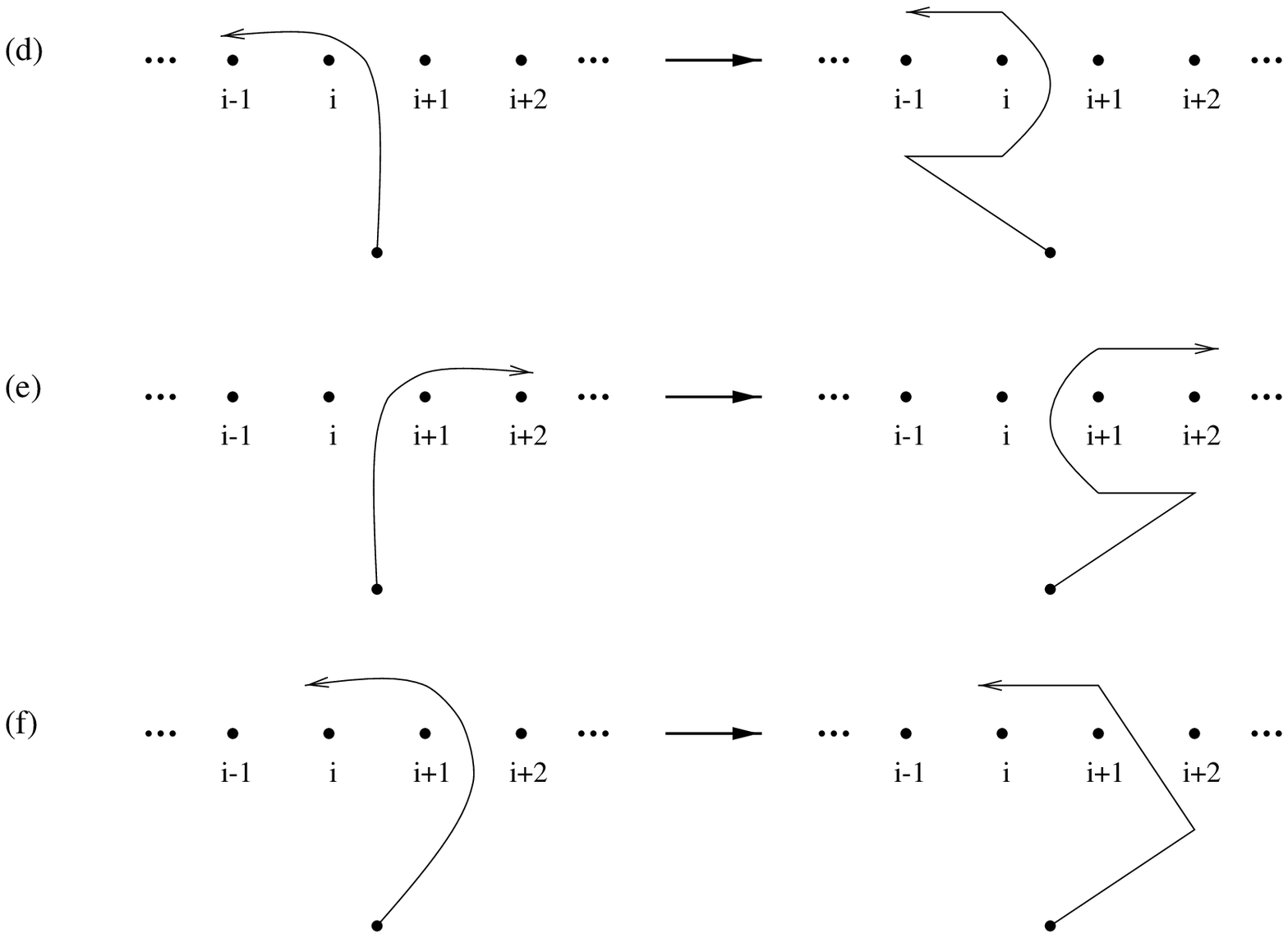}

\small{Figure $12$: Homotopical modifications of the elements of the g-base}
\end{center}

This concludes the proof of the correctness of the $PositiveHalfTwist(\sigma _i,g)$ and the 
$NegativeHalfTwist(\sigma _i,g)$ functions. 
We still have to prove the correctness of the $Reduce(g)$ function, and
that it does not change the homotopy type of the elements of the g-base. These proofs will make it possible 
to derive the uniqueness of the presentation.

\subsection{Correctness of the $Reduce(g)$ algorithm}
Here we will lay out the proof of the correctness of the $Reduce(g)$ algorithm. We will prove that 
the algorithm does not change the homotopy type of the elements of the g-base, and that it returns a list which 
represents the g-base in a unique form.

\bpr
Let $g$ be a list representing a g-base. Then, the list returned by the function
$Reduce(g)$ represents a g-base which is homotopically equivalent to $g$.
\epr

\noindent
{\it Proof:} The $Reduce(g)$ algorithm is based on four reduction rules. We will present
the rules and we will prove that each one of them preserves the homotopy type of $g$.

\ben
\item
If we have two consecutive equal links, we can omit them both.
\item
If we have a sequence of $(i,\pm 1) \to (i,0)$, we can omit the first link.
\item
If we have links between $(i,0)$ and $(-1,0)$, we can omit them all.
\item
If we have a sequence that starts with $(-1,0)$ and continues to $(i,-1)$, we can omit the
latter link.
\een

Concerning the first rule, the meaning of the situation of two consecutive equal links, is
that the path is moving above (or beneath) a point and immediately retracing back. 
Homotopically, this is equivalent to a point. Hence, we can omit the two links.

Concerning the second rule, the link $(i,\pm 1)$ represents a point which is directly
above or below the point $(i,0)$ and very close to it. Therefore, we result in a homotopic path after omitting
the link $(i,\pm 1)$.

The third rule is trivial since any link that was added between the end point of one path
$(i,0)$ and the beginning point of the next path $(-1,0)$ is not even a part of the g-base 
presentation and therefore has to be erased.

The fourth rule is based on the fact that the shortest path between two points is a straight line. 
Therefore, any point that we add can be omitted
without changing the homotopy type of the g-base. We should remember that in our
presentation of the paths of the g-base, we use the convention that the start point
$(-1,0)$ will always be connected directly to the point $(i,0)$ or to the point above
$(i,1)$ but never to the point $(i,-1)$.
\hfill $\qed$

\ble
Let $g$ be a representation of a g-base. Then $g$ does not contain any sequence of the
following type $(i-1,e) \to (i,\pm 1) \to (i,\mp 1) \to (i+1,e)$, $e \in \{-1,1,0\}$.
\ele

\noindent
{\it Proof:} We will prove it by induction. The initial g-base we have is the
standard g-base, which does not contain any such sequence.

Now, at each step, when we add new links we connect them to links to the left or right of the local
twisted section. This means that we will always
connect the local section of the path to the other ends by a sequence that does not
contain both $(i,-1)$ and $(i,1)$, and therefore we will not create at any step a sequence
that was forbidden by the lemma. Consider the fact that the change in the local section of the path 
is only a twist of $180 ^\circ$, then this twist will not add any forbidden sequences.
Moreover, using the $Reduce(g)$ function, we eliminate any
unnecessary links, resulting in the fact that each connection is the shortest possible (with regard to
the convention that we will never connect the point $(-1,0)$ to any $(i,-1)$). Therefore, the
$Reduce(g)$ function will not add any forbidden sequences. So by induction we proved the
lemma.
\hfill $\qed$

\bco
Let $g$ be a representation of a g-base. Then the representation of the g-base obtained by 
$Reduce(g)$ is unique.
\eco

\noindent
{\it Proof:} As stated above, by the correctness of the algorithms, and since $Reduce(g)$
will always connect two points by the shortest path without changing the homotopic type of
the g-base, we obtain a unique representation of the g-base.
\hfill $\qed$

This concludes the proof of the correctness of the algorithm.

\bth
The $ProcessWord(w)$ algorithm will result in a unique representation of the g-base 
after the action of the braid word $w$ on the standard g-base.
\eth
\hfill $\qed$

\section{Complexity}
In this section, we will compute the complexity of the two functions that we use: 
$(Positive/Negative)HalfTwist(\sigma _i,g)$ and $Reduce(g)$.

\bpr
Let $l_g$ be the length of the list representing a g-base. Then, the complexity of the function
$(Positive/Negative)HalfTwist(\sigma _i,g)$ is bounded by $O(l_g)$.
\epr

\noindent
{\it Proof:} The algorithm goes through all the links in the list, acting at most once on every link.
As an upper bound, the algorithm might add two links for every link that was in the
list, resulting in a list with length $3l_g$. Therefore at most the algorithm will perform
$4l_g$ operations, which yields in the complexity of $O(l_g)$.
\hfill $\qed$

\brem
In practice the number of links actually added is much smaller than the upper bound given above. 
This is the reason we have a practical and an efficient algorithm for short words.
\erem

\bpr
Let $l_g$ be the length of the list representing a g-base. Then, the complexity of the function
$Reduce(g)$ is bounded by $O(l_g)$.
\epr

\noindent
{\it Proof:} This result is a consequence of the fact that the algorithm will check every link at
most twice, and that each link that was inserted in the list can be extracted only once. We
have to notice that there cannot be a situation of the following type $(i,e) \to (i+1,e) \to (i+2,e)
\to (i+2,e) \to (i+1,e) \to (i,e)$, since after each step of inserting links we delete the
unnecessary ones, and the fact that the local section is always two points wide. 
This fact is what makes it possible to make sure that in the worst case while
going through the list, in order to delete unnecessary links, we will need to retrace only one step. 
That means that keeping in memory the link before the one we are currently checking is sufficient, and
that each link can be passed at most twice; therefore, the complexity bound is $O(l_g)$.
\hfill $\qed$

Note that even if one wants to use a general half-twist which allows the local 
section to be larger than two points (with proper modification of the prefix and 
postfix sequences), obtaining the complexity of $O(l_g)$ is still possible by using 
a doubly-connected list.

\section{Conclusions}
We would like to state here that although for very long braid words this algorithm's running time is long, 
since the complexity of the g-base presentation grows with the length of the braid word, for short braid words
we have obtained a quick algorithm in comparison with other methods. This is true because of
the fact that Garside's algorithm involves the replacement of the generators in a negative
power by a subword of size $n(n-1)$, where $n$ is the number of strings in the
braid group (although, for variations of his algorithm the size of the fundamental word $\Delta _n$ reduces \cite{NEW}), 
and because each step of reduction in Dehornoy's algorithm involves the
insertion of at least two subwords of $O(n)$ length. In our algorithm for short braid
words over a large number of strings we obtain a very short description of the g-base, that
yields a very fast algorithm. This means that we have presented a very useful and practical
algorithm. We also would like to point out that we have an implementation of the algorithm 
on a computer.

Instead of the restriction on the length of the braid word, one might consider a
restriction on the size of the presentation of the g-base. Therefore by excluding braid
word where the number of twists of the g-base's paths grows dramatically, we still have an
efficient algorithm even for longer braid words.

We would like to point out some of the future applications in which we believe that the new
approach may help. We think that there has to be a connection between presentations 
of two conjugated braid words, what might bring a practical fast algorithm for solving the 
conjugacy problem in the braid group. 

We believe that we can do the unprocess of the algorithm, which means to compute the braid word from a given g-base.
Another thing is obvious when looking at the braid monodromy as a homomorphism from one
fundamental group of a punctured disk to another. We believe that using this new method 
will make it easier to compute the braid monodromy automatically, at least for
some of the cases.

Another implication of the algorithm and the new method for presenting the g-base is a
similar implementation for the Moishezon-Teicher algorithm for computing the braid
monodromy of real line arrangements and plane curves. 
For this implementation we have an efficient computer program.

\begin{\bib}{10}
\bibitem{Artin} E. Artin, Theory of braids, {\it Ann. Math.} {\bf 48} (1947), 101-126.
\bibitem{NEW} J.S. Birman, K.H. Ko and S.J. Lee, A new approach to the word and conjugacy problems in the braid groups, {\it Adv. Math.} {\bf 139} (1998), 322-353.
\bibitem{Deh1} P. Dehornoy, From large cardinals to braids via distributive algebra, {\it J. Knot Theory \& Ramifications} {\bf 4(1)} (1995), 33-79.
\bibitem{Deh2} P. Dehornoy, A fast method for comparing braids, {\it Adv. Math.} {\bf 125(2)} (1997), 200--235. 
\bibitem{POSBR} E.A. Elrifai and H.R. Morton, Algorithms for positive braids, {\it Quart. J. Math. Oxford Ser. (2)} {\bf 145} (1994), 479-497.
\bibitem{GAR} F.A. Garside, The braid group and other groups, {\it Quart. J. Math. Oxford Ser. (2)} {\bf 78} (1969), 235-254.
\bibitem{EFF} A. Jacquemard, About the effective classification of conjugacy classes of braids, {\it J. Pure. Appl. Alg.} {\bf 63} (1990), 161-169.
\bibitem{BAND} E.S. Kang, K.H. Ko and S.J. Lee, Band-generator presentation for the 4-braid group, {\it Top. Appl.} {\bf 78} (1997), 39-60.
\bibitem{BGTI} B. Moishezon and M. Teicher, Braid group techniques in complex geometry I, Line arrangements in $\C \PP^2$, {\it Contemporary Math.} {\bf 78} (1988), 425-555.
\end{\bib}

\end{document}